\documentclass[11pt,a4paper]{amsart}
\usepackage{amssymb,amsmath}
\usepackage{mathpazo}
\headsep=1cm \numberwithin{equation}{section}
\newtheorem{theorem}{Theorem}[section]
\newtheorem{lemma}[theorem]{Lemma}

\newtheorem{corollary}[theorem]{Corollary}

\theoremstyle{definition}
\newtheorem{definition}[theorem]{Definition}
\theoremstyle{remark}
\newtheorem{remark}[theorem]{Remark}
\newtheorem{example}[theorem]{Example}

\newcommand{\Ext}{\operatorname{Ext}}

\newcommand{\Tor}{\operatorname{Tor}}
\newcommand{\Hom}{\operatorname{Hom}}

\newcommand{\Ann}{\operatorname{Ann}}

\newcommand{\lo}{\longrightarrow}
\newcommand{\fm}{\frak{m}}
\newcommand{\fp}{\frak{p}}

\newcommand{\fa}{\frak{a}}

\def\mapdown#1{\Big\downarrow\rlap
{$\vcenter{\hbox{$\scriptstyle#1$}}$}}

\begin{document}
\author[Divaani-Aazar, Esmkhani and Tousi ]{Kamran Divaani-Aazar,
Mohammad Ali Esmkhani and Massoud Tousi}
\title[A criterion for rings which are locally valuation]
{A criterion for rings which are locally valuation}

\address{K. Divaani-Aazar, Department of Mathematics, Az-Zahra
University, Vanak, Post Code 19834, Tehran, Iran-and-Institute for
Studies in Theoretical Physics and Mathematics, P.O. Box
19395-5746, Tehran, Iran.} \email{kdivaani@ipm.ir}

\address{M.A. Esmkhani, Department of Mathematics, Shahid Beheshti
University, Tehran, Iran-and-Institute for Studies in Theoretical
Physics and Mathematics, P.O. Box 19395-5746, Tehran, Iran.}

\address{M. Tousi, Department of Mathematics, Shahid Beheshti
University, Tehran, Iran-and-Institute for Studies in Theoretical
Physics and Mathematics, P.O. Box 19395-5746, Tehran, Iran.}


\keywords{Absolutely pure modules, cyclically pure injective
modules, projective principal rings, Pr\"{u}fer domains,
semi-hereditary rings, semi-simple rings,
valuation rings.\\
The first author was supported by a grant from IPM (No. 84130213).\\
The third author was supported by a grant from IPM (No.
85130213).}

\begin{abstract} Using the notion of cyclically pure
injective modules,  a characterization for rings which are locally
valuation is established. As applications, new characterizations
for Pr\"{u}fer domains and pure semi-simple rings are provided.
Namely, we show that a domain $R$ is Pr\"{u}fer if and only if two
of the three classes of pure injective, cyclically pure injective
and RD-injective modules are equal. Also, we prove that a
commutative ring $R$ is pure semi-simple if and only if every
$R$-module is cyclically pure injective.
\end{abstract}

\maketitle

\section{Introduction}

Throughout this paper, R denotes a commutative ring with identity,
and all modules are assumed to be left unitary. The notion of pure
injective modules has a substantial role in commutative algebra
and model theory. Among various generalizations of this notion,
the notion of cyclically pure injective modules has been
extensively studied by M. Hochster [{\bf 9}] and L. Melkersson
[{\bf 12}]. Recall that an exact sequence $0\lo A\lo B\lo C\lo 0$
of $R$-modules and $R$-homomorphisms is said to be {\it cyclically
pure} if the induced map $R/\fa\otimes_RA\lo R/\fa\otimes_RB$ is
injective for all (finitely generated) ideals $\fa$ of $R$. Also,
an $R$-module $D$ is said to be {\it cyclically pure injective} if
for any cyclically pure exact sequence $0\lo A\lo B\lo C\lo 0,$
the induced homomorphism $\Hom_R(B,D)\lo \Hom_R(A,D)$ is
surjective. In the sequel, we use the abbreviation CP for the term
``cyclically pure".

More generally, let $\mathcal{S}$ be a class of $R$-modules. An
exact sequence $0\lo A\lo B\lo C\lo 0$ of $R$-modules and
$R$-homomorphisms is said to be $\mathcal{S}$-pure if for all $M\in
\mathcal{S}$, the induced homomorphism $\Hom_R(M,B)\lo \Hom_R(M,C)$
is surjective. An $R$-monomorphism $f:A\lo B$ is said to be
$\mathcal{S}$-pure if the exact sequence $0\lo A\overset{f}\lo
B\overset{nat}\lo B/f(A)\lo 0$ is $\mathcal{S}$-pure. An $R$-module
$D$ is said to be $\mathcal{S}$-pure injective if for any
$\mathcal{S}$-pure exact sequence $0\lo A\lo B\lo C\lo 0,$ the
induced homomorphism $\Hom_R(B,D)\lo \Hom_R(A,D)$ is surjective, see
[{\bf 13}]. When $\mathcal{S}$ is the class of finitely presented
$R$-modules, $\mathcal{S}$-pure exact sequences and
$\mathcal{S}$-pure injective modules are called pure exact sequences
and pure injective modules, respectively. If $\mathcal{S}$ denotes
the class of all $R$-modules of the form $R/Rr,r\in R$, then
$\mathcal{S}$-pure exact sequences and $\mathcal{S}$-pure injective
modules are called RD-exact sequences and RD-injective modules,
respectively. For a survey on the notions of pure injective and
RD-injective modules, we refer the reader to [{\bf 6}].

Let $\mathcal{S}$ be the class of all $R$-modules $M$ for which
there is a cyclic submodule $G$ of $R^{n}$, for some $n\in
\mathbb{N}$, such that $M$ is isomorphic to $R^n/G$. In [{\bf 3}],
we showed that CP-exact sequences and CP-injective modules coincide
with $\mathcal{S}$-pure exact sequences and $\mathcal{S}$-pure
injective modules, respectively. In the same paper we have
systematically investigated the structure of CP-injective modules
and presented several characterizations of this class of modules.
Our aim in this paper is the following:

i) Classifying the commutative rings that over which the two notions
of ``RD-injective" and ``cyclically pure injective" coincide.

ii) Classifying the commutative rings that over which the two
notions of ``pure injective" and ``cyclically pure injective"
coincide.

In Section 2, we show that $R_{\fp}$ is a valuation ring for all
prime ideals $\fp$ of $R$ if and only if every CP-injective
$R$-module is RD-injective, if and only if every pure injective
$R$-module is CP-injective. From this we obtain a characterization
for semi-hereditary rings and also one for Pr\"{u}fer domains. In
the literature, there are several characterizations for Pr\"{u}fer
domains. In particular, by [{\bf 6}, Chapter XIII, Theorem 2.8], it
is known that a domain $R$ is Pr\"{u}fer if and only if every pure
injective $R$-module is RD-injective. Also, it is known by [{\bf 6},
Chapter IX, Proposition 3.4] that a domain $R$ is Pr\"{u}fer if and
only if every divisible $R$-module is absolutely pure. Here we show
that a domain $R$ is Pr\"{u}fer if and only if every CP-injective
$R$-module is RD-injective, if and only if every pure injective
$R$-module is CP-injective.  Also, we show that a domain $R$ is
Pr\"{u}fer if and only if every absolutely CP-module is absolutely
pure. Finally, a new characterization for pure semi-simple rings is
given. We show that a ring $R$ is pure semi-simple if and only if
every $R$-module is CP-injective, if and only if every $R$-module is
RD-pure injective.

The first example of a CP-exact sequence which is not pure was
presented in [{\bf 1}]. Our first characterization for Pr\"{u}fer
domains mentioned above shows that over a non-Pr\"{u}fer domain
$R$ the class of CP-injective $R$-modules is strictly larger than
that of RD-injective $R$-modules and strictly smaller than that of
pure injective $R$-modules. However, these may be viewed as kind
of implicit strict inclusions. In Section 3, we provide some
examples for which we can explicitly show proper containments in
this regard. In [{\bf 3}], we proved that in many aspects
CP-injective modules behave similar to pure injective and
RD-injective modules. But Remark 2.2 and Example 3.5 below display
some differences between the former class and the later two.

\section{A characterization for Pr\"{u}fer rings}

In the remainder of this paper, let $\mathcal{S}_1$ denote the
class of all $R$-modules of the form $R/Rr,r\in R$. Also, let
$\mathcal{S}_4$ (resp. $\mathcal{S}_2$) denote the class of all
finitely presented (resp. finitely presented cyclic) $R$-modules.
Finally, we let $\mathcal{S}_3$ denote the class of all
$R$-modules $M$ for which there are an integer $n\in \mathbb{N}$
and a cyclic submodule $G$ of $R^ n$ such that $M$ is isomorphic
to $R^ n/G$.

\begin{definition} Let $\mathcal{S}$ be a class of $R$-modules.
An exact sequence $0\lo A\lo B\lo C\lo 0$ of $R$-modules and
$R$-homomorphisms is called $\mathcal{S}$-flat if for all $M\in
\mathcal{S}$ the induced map $A\otimes_RM\lo B\otimes_RM$ is
injective.
\end{definition}

\begin{remark} Let $0\lo A\lo B\lo C\lo 0$ be an exact sequence of
$R$-modules and $R$-homomorphisms.\\
i) For $i=1, 4$, the above exact sequence is $\mathcal{S}_i$-pure
if and only if it is $\mathcal{S}_i$-flat, see [{\bf 13},
Propositions 2 and 3].\\
ii) By [{\bf 3}, Proposition 2.2], the above exact sequence is
$\mathcal{S}_3$-pure if and only if it is $\mathcal{S}_2$-flat.
\end{remark}

Example 3.5 in the next section, shows that there exists an
$\mathcal{S}_2$-flat exact sequence which is not
$\mathcal{S}_2$-pure.

\begin{definition} Let $\mathcal{S}$ be a class of $R$-modules.
An $R$-module $P$ is said to be $\mathcal{S}$-{\it pure
projective} if for any $\mathcal{S}$-pure exact sequence $0\lo
A\lo B\lo C\lo 0 ,$ the induced homomorphism $\Hom_R(P,B)\lo
\Hom_R(P,C)$ is surjective.
\end{definition}

\begin{lemma} Let $\mathcal{S}$ and $\mathcal{T}$ be two classes
of $R$-modules. The following are equivalent:\\
i) Every $\mathcal{T}$-pure exact sequence is $\mathcal{S}$-pure
exact.\\
ii) Every $\mathcal{S}$-pure projective $R$-module is
$\mathcal{T}$-pure projective.\\
iii) Every element of $\mathcal{S}$ is a direct summand of a
direct sum of modules in $\mathcal{T}$.\\
Moreover, if $\mathcal{S}$ and $\mathcal{T}$ are both contained in
$\mathcal{S}_4$, then the above conditions  are equivalent to the
following\\
iv) Every $\mathcal{S}$-pure injective $R$-module is
$\mathcal{T}$-pure injective.
\end{lemma}

{\bf Proof.} Let $\mathcal{U}$ be a class of $R$-modules. By the
definition every element of $\mathcal{U}$ is $\mathcal{U}$-pure
projective. In general, by [{\bf 13}, Proposition 1], it turns out
that an $R$-module $M$ is $\mathcal{U}$-pure projective if and
only if $M$ is a direct summand of a direct sum of modules in
$\mathcal{U}$. Hence the equivalence of i), ii) and iii) is
immediate.

Next, assume that $\mathcal{S}$ and $\mathcal{T}$ are both
contained in $\mathcal{S}_4$. Let $\mathcal{U}\subseteq
\mathcal{S}_4$ be a class of $R$-modules and $E$ an injective
cogenerator of $R$. By [{\bf 5}, Lemma 1.2], there is a class
$\mathcal{U}^*$ of $R$-modules such that an exact sequence $0\lo
A\lo B\lo C\lo 0$ of $R$-modules and $R$-homomorphisms is
$\mathcal{U}$-pure if and only if $$0\lo A\otimes_RM^*\lo
B\otimes_RM^*\lo C\otimes_RM^*\lo 0$$ is exact for all $M^*\in
\mathcal{U}^*$. Thus by using adjoint property, it follows that
$\Hom_R(M^*,E)$ is a $\mathcal{U}$-pure injective $R$-module for
all $M^*\in \mathcal{U}^*$.

$iv)\Rightarrow i)$ Let $0\lo A\lo B\lo C\lo 0 (*)$ be a
$\mathcal{T}$-pure exact sequence and $M^*\in \mathcal{S}^*$ an
arbitrary element. Since $\Hom_R(M^*,E)$ is $\mathcal{S}$-pure
injective, it is also  $\mathcal{T}$-pure injective, by our
assumption. Thus, by applying the functor
$\Hom_R(-,\Hom_R(M^*,E))$ on $(*)$ and using adjoint property, we
deduce the following exact sequence
$$ 0\lo \Hom_R(C\otimes_RM^*,E)\lo \Hom_R(B\otimes_RM^*,E)\lo
\Hom_R(A\otimes_RM^*,E)\lo 0.$$ Thus, it turns out that the
sequence
$$0\lo A\otimes_RM^*\lo B\otimes_RM^*\lo C\otimes_RM^*\lo 0$$ is
exact. Therefore $(*)$ is $\mathcal{S}$-pure exact.

Now, since the implication $i)\Rightarrow iv)$ clearly holds, the
proof is finished. $\Box$

\begin{lemma} Assume that every pure injective $R$-module is
CP-injective. Then an exact sequence $l: 0\lo A\lo B\lo C\lo 0$ is
$\mathcal{S}_2$-pure exact if and only if it is CP-exact.
\end{lemma}

{\bf Proof.}  Assume that $l$ is a CP-exact sequence. Then, by
Lemma 2.4, it is pure exact. Hence it is clearly
$\mathcal{S}_2$-pure, because
$\mathcal{S}_2\subseteq\mathcal{S}_4$.

Now, assume that $l$ is $\mathcal{S}_2$-pure exact. Let $E$ be an
injective cogenerator of $R$ and $(\cdot)^ \vee$ denote the
faithfully exact functor $\Hom_R(-,E)$. Let $l^ \vee$ denote the
induced exact sequence $0\lo C^ \vee\lo B^ \vee\lo A^ \vee\lo 0$.
Let $\frak{I}$ be a finitely generated ideal of $R$. Since
$R/{\frak{I}}$ is finitely presented, the two $R$-modules
$R/{\frak{I}}\otimes_R{M^ \vee}$ and ${\Hom_R(R/{\frak{I}},M)}^
\vee$ are naturally isomorphic for all $R$-modules $M$. So the
exact sequence $l^ \vee$ is a CP-exact. Hence $l^ \vee$ is pure
exact, by Lemma 2.4. Let $N\in \mathcal{S}_3$. Then by Remark 2.2
i), the sequence $N\otimes_Rl^ \vee$ is exact. The exact sequences
$$0\lo N\otimes_RC^ \vee\lo N\otimes_RB^ \vee\lo N\otimes_RA^
\vee\lo 0$$ and $$0\lo {\Hom_R(N,C)}^ \vee \lo {\Hom_R(N,B)}^ \vee
\lo {\Hom_R(N,A)}^ \vee \lo 0$$ are naturally isomorphic. Thus the
second sequence is also exact, and so
$$0\lo \Hom_R(N,A)\lo \Hom_R(N,B)\lo \Hom_R(N,C)\lo 0$$ is an
exact sequence, because $(\cdot)^ \vee$ is a faithfully exact
functor. Therefore $l$ is a CP-exact sequence. $\Box$

\begin{lemma} Let $\fa$ be an ideal of $R$. Assume that every
CP-injective $R$-module is RD-injective. Then every CP-injective
$R/{\fa}$-module is an RD-injective $R/{\fa}$-module.
\end{lemma}

{\bf Proof.} Set $T=R/{\fa}$. Let $M={T^ n}/V$, where $n\in
\mathbb{N}$ and $V$ is a cyclic $T$-submodule of ${T^ n}$. So,
there are $b_1, \ldots, b_n \in R$ such that $V=T(b_1+\fa, \ldots,
b_n+\fa)$. Let $N={R^ n}/U$, where $U=R(b_1, \ldots, b_n)$. We
show that $M$ and $N\otimes_RT$ are naturally isomorphic as
$T$-modules. To this end, let $\phi:M\lo N\otimes_RT$ be the map
defined by
$$(x_1+\fa, \ldots, x_n+\fa)+V \mapsto ((x_1, \ldots,
x_n)+U)\otimes(1+\fa)
$$ for all $(x_1+\fa,\ldots, x_n+\fa)+V \in M$. Also, we define
$\psi:N\otimes_RT\lo M$ by $$((x_1, \ldots,
x_n)+U)\otimes(r+\fa)\mapsto (rx_1+\fa, \ldots, rx_n+\fa)+V.
$$ It is a routine check to see that $\phi$ and $\psi$ are well defined
$T$-homomorphisms and that $\psi \phi=id_M$ and $\phi
\psi=id_{N\otimes_RT}$. Now, as $-\otimes_RT$ commutes with direct
sums, the conclusion is immediate by Lemma 2.4
$iii)\Longleftrightarrow iv)$. $\Box$

Recall that a {\it valuation ring} (not necessarily a domain) is a
commutative ring whose ideals are linearly ordered under
inclusion.

\begin{theorem} The following are equivalent:\\
i) $R_{\fp}$ is a valuation ring for all prime ideals $\fp$ of $R$.\\
ii) Every pure injective $R$-module is RD-injective.\\
iii) Every CP-injective $R$-module is RD-injective.\\
iv) Every pure injective $R$-module is CP-injective.\\
v) Every pure projective $R$-module is RD-projective.\\
vi) Every CP-projective $R$-module is RD-projective.\\
vii) Every pure projective $R$-module is CP-projective.
\end{theorem}

{\bf Proof.} By Lemma 2.4, the equivalences
$ii)\Longleftrightarrow v)$, $iii)\Longleftrightarrow vi)$ and
$iv)\Longleftrightarrow vii)$ are obvious. Also, the implications
$ii)\Rightarrow iii)$ and $ii)\Rightarrow iv)$ are clear.

$i)\Rightarrow v)$ As we have mentioned in the proof Lemma 2.4,
for a given class $\mathcal{U}$ of $R$-modules,  an $R$-module $M$
is $\mathcal{U}$-pure projective if and only if $M$ is a direct
summand of a direct sum of modules in $\mathcal{U}$. So, to deduce
v), it is enough to show that every finitely presented $R$-module
is RD-projective. By [{\bf 6}, Proposition 4], a finitely
presented $R$-module $M$ is RD-projective if and only if $M_{\fm}$
is an RD-projective $R_{\fm}$-module for all maximal ideals $\fm$
of $R$. Hence v) follows by [{\bf 15}, Theorem 1].

$v)\Rightarrow i)$ follows by [{\bf 13}, Proposition 1] and [{\bf
15}, Theorem 3].

$iii)\Rightarrow i)$ Assume that there exists a prime ideal $\fp$
of $R$ so that $R_{\fp}$ is not a valuation ring. Let
$N=(R_{\fp})^n/G$, where $n\in \mathbb{N}$ and $G$ is a cyclic
$R_{\fp}$-submodule of $(R_{\fp})^n$. Clearly $N$ is equal to the
localization at $\fp$ of an element of $\mathcal{S}_3$. Hence, as
localization at $\fp$ commutes with direct sums, by Lemma 2.4, we
may and do assume that $R$ is a local ring which is not a
valuation ring. Denote by $\fm$ the maximal ideal of $R$. Since
$R$ is not a valuation ring, there are two elements $a,b\in R$
such that $Ra\nsubseteq Rb$ and $Rb\nsubseteq Ra$. Set
$\frak{I}:=\fm a+\fm b$. Lemma 2.6 yields that every CP-injective
$R/{\frak{I}}$-module is an RD-injective $R/{\frak{I}}$-module.
Replace $R, a$ and $b$ by $R/{\frak{I}}, a+{\frak{I}}$ and
$b+{\frak{I}}$, respectively. So we can assume that $R$ is a local
ring which is not a valuation ring and that there are two elements
$a,b\in R$ such that $Ra\nsubseteq Rb$, $Rb\nsubseteq Ra$, $\fm
a=\fm b=0$ and $Ra\cap Rb=0$. In view of the proof of [{\bf 15},
Theorem 2], it becomes clear that $M:=(R\oplus R)/R(a, -b)$ is a
non-cyclic indecomposable $R$-module. Lemma 2.4 implies that $M$
is a direct summand of a direct sum of cyclic modules. Now, by
[{\bf 14}, Proposition 3], over a commutative local ring, any
indecomposable direct summand of a direct sum of cyclic modules is
cyclic. We achieved at a contradiction.

$iv)\Rightarrow i)$ By Lemmas 2.4 and 2.5, it follows that every
finitely presented $R$-module is a direct summand of a direct sum
of cyclic modules. Now, we assume that i) does not hold and search
for a contradiction. Then there is a prime ideal $\fp$ of $R$ so
that $R_{\fp}$ is not a valuation ring. Hence, by [{\bf 15},
Theorem 2], there exists an indecomposable finitely presented
$R_{\fp}$-module $M$ which is not cyclic. Since every finitely
presented $R_{\fp}$-module is the localization at $\fp$ of a
finitely presented $R$-module, we deduce that $M$ is a direct
summand of a direct sum of cyclic $R_{\fp}$-module. But then by
[{\bf 14}, Proposition 3], $M$ should be a cyclic
$R_{\fp}$-module. $\Box$

\begin{definition} i) (See [{\bf 4}]) A ring $R$ is said to be {\it
projective principal} ring (P.P.R.) if every principal ideal of $R$
is projective.\\
ii) A ring $R$ is said to be {\it semi-hereditary} if every
finitely generated ideal of $R$ is projective.\\
iii) (See [{\bf 10}]) An $R$-module $M$ is said to be {\it
absolutely pure} (resp. {\it absolutely cyclically pure}) if it is
pure (resp. cyclically pure)
as a submodule in every extension of $M$.\\
iv) (See [{\bf 5}]) An $R$-module $M$ is said to be {\it
divisible} if for every $r\in R$ and $x\in M$, $\Ann_Rr\subseteq
\Ann_Rx$ implies that $x\in rM$. (This is equivalent to the usual
definition where $R$ is domain.)
\end{definition}

In the proof of the following lemma we use the methods of the
proofs of [{\bf 10}, Proposition 1 and Corollary 2].

\begin{lemma} Let $M$ be an $R$-module.\\
i) $M$ is absolutely cyclically pure if and only if
$\Ext_R^ i(N,M)=0$ for all $N\in \mathcal{S}_3$.\\
ii) $M$ is absolutely cyclically pure if and only if any diagram
\begin{equation*}
\setcounter{MaxMatrixCols}{11}
\begin{matrix}
&P' &\overset{\alpha} \lo  &P
\\&\mapdown{\beta} & &
& & & &
\\&M
\end{matrix}
\end{equation*}
with $P'$ cyclic, $\alpha$ monic and $P$ projective, there exists
a homomorphism $\gamma:P\lo M$ such that $\gamma \alpha= \beta$.
\end{lemma}

{\bf Proof.} i) Let $L$ be an extension of $M$ and $N\in
\mathcal{S}_3$. From the exact sequence $0\lo M\hookrightarrow
L\lo L/M\lo 0$, we deduce the following exact sequence
\begin{small}
$$0\rightarrow \Hom_R(N,M)\rightarrow \Hom_R(N,L)\rightarrow
\Hom_R(N,L/M)\rightarrow \Ext_R^1(N,M)\rightarrow \Ext_R^1(N,L)
(*).$$
\end{small}
Assume that $M$ is an absolutely CP-module and
let $L$ be an injective extension of $M$. Then by Remark 2.2 ii)
and $(*)$, we conclude that $\Ext_R^1(N,M)=0$ for all $N\in
\mathcal{S}_3$.

Now, assume that $\Ext_R^1(N,M)=0$ for all $N\in \mathcal{S}_3$.
Let $L$ be an extension of $M$. Then Remark 2.2 ii) and $(*)$
imply that the exact sequence $0\lo M\hookrightarrow L\lo L/M\lo
0$ is CP-exact.

ii) We may assume that $P$ is a finitely generated free
$R$-module. Thus the result follows by using i) and the following
exact sequence $$\Hom_R(P,M)\lo \Hom_R(P',M)\lo \Ext_R^
1(P/\alpha(P'),M)\lo 0.\Box$$

\begin{lemma} The following are equivalent:\\
i) $R$ is a P.P.R.\\
ii) Every cyclic submodule of a projective $R$-module is
projective.\\
iii) Every quotient of an absolutely CP-module is also an
absolutely CP-module.
\end{lemma}

{\bf Proof.} $i)\Leftrightarrow ii)$ follows by [{\bf 4}, Theorem
3.2].

$ii)\Leftrightarrow iii)$ In view of Lemma 2.9, the proof is
immediate  by adapting  the argument of [{\bf 10}, Theorem 2] and
replacing the phrases ``absolutely pure" and ``finitely generated
submodule" with ``absolutely cyclically pure" and ``cyclic
submodule", respectively. $\Box$

\begin{corollary} Assume that $R$ is a P.P.R. The following are
equivalent:\\
i) $R$ is a semi-hereditary ring.\\
ii) Every pure injective $R$-module is RD-injective.\\
iii) Every CP-injective $R$-module is RD-injective.\\
iv) Every pure injective $R$-module is CP-injective.\\
v) Every divisible $R$-module is absolutely pure.\\
vi) Every absolutely CP-module is absolutely pure.\\
vii) Every pure projective $R$-module is RD-projective.\\
viii) Every CP-projective $R$-module is RD-projective.\\
ix) Every pure projective $R$-module is CP-projective.
\end{corollary}

{\bf Proof.} As, we have mentioned in the proof of Theorem 2.7, by
Lemma 2.4, the equivalences $ii)\Longleftrightarrow vii)$,
$iii)\Longleftrightarrow viii)$ and $iv)\Longleftrightarrow ix)$
are obvious.

Now, assume that $R$ is semi-hereditary. Let $\fp$ be a prime
ideal of $R$. Then $R_{\fp}$ is also a semi-hereditary ring. Hence
for each nonzero element $a$ of $R_{\fp}$, the $R_{\fp}$-module $a
R_{\fp}$ is a nonzero free $R_{\fp}$-module. Thus, we conclude
that $R_{\fp}$ is a domain. But, it is known that a domain is
semi-hereditary if and only if it is Pr\"{u}fer. So $R_{\fp}$ is a
valuation domain for all prime ideals $\fp$ of $R$. Therefore the
implication $i)\Rightarrow ii)$ and the equivalences
$ii)\Leftrightarrow iii)$ and $iii)\Leftrightarrow iv)$ are
immediate by Theorem 2.7.

$ii)\Rightarrow v)$ Let $M$ be a divisible $R$-module and $E$
denote the injective envelop of $M$. Then [{\bf 5}, Lemma 2.2]
implies that the sequence $0\lo M\hookrightarrow E\lo E/M\lo 0$,
is RD-exact. Hence, by Lemma 2.4, it is pure and so $\Ext_R^
1(N,M)=0$ for all $N\in \mathcal{S}_4$. Thus, by [{\bf 10},
Proposition 1], $M$ is absolutely pure.

$v)\Rightarrow vi)$ Let $M$ be an absolutely CP-module. Then, by
Lemma 2.9 i), $\Ext_R^ 1(N,M)=0$ for all $N\in \mathcal{S}_3$. In
particular, $\Ext_R^ 1(R/Rr,M)=0$ for all $r\in R$, and so $M$ is
a divisible $R$-module by [{\bf 5}, Lemma 2.2]. Thus $M$ is
absolutely pure, as required.

Finally, we prove $vi)\Rightarrow i)$. Since $R$ is a P.P.R.,
Lemma 2.10 yields that every quotient of an absolutely CP-module
is again an absolutely CP-module. So, if vi) holds, then every
quotient of an absolutely pure module is again absolutely pure.
Thus i) follows by [{\bf 10}, Theorem 2]. $\Box$

Now, since a domain $R$ is Pr\"{u}fer if and only if it is
semi-hereditary, we can obtain the main result of this paper. Note
that every domain is a P.P.R.

\begin{corollary} Assume that $R$ is a domain. The following are
equivalent:\\
i) $R$ is Pr\"{u}fer.\\
ii) Every pure injective $R$-module is RD-injective.\\
iii) Every CP-injective $R$-module is RD-injective.\\
iv) Every pure injective $R$-module is CP-injective.\\
v) Every divisible $R$-module is absolutely pure.\\
vi) Every absolutely CP-module is absolutely pure.\\
vii) Every pure projective $R$-module is RD-projective.\\
viii) Every CP-projective $R$-module is RD-projective.\\
ix) Every pure projective $R$-module is CP-projective.
\end{corollary}

Let $\mathcal{C}_{RDR}$ denote the class of all RD-injective
$R$-modules. Also, let $\mathcal{C}_{CPR}$ and $\mathcal{C}_{PR}$
denote the class of all CP-injective $R$-modules and that of all
pure injective $R$-modules, respectively. It follows, by Theorem
2.7 that if two of three classes $\mathcal{C}_{RDR}$,
$\mathcal{C}_{CPR}$ and $\mathcal{C}_{PR}$ are equal, then all
three classes are equal. The following result shows that if each
of these three classes is equal to the class of all $R$-modules,
then the two other classes are also equal to the class of all
$R$-modules. First, we bring a definition.

\begin{definition} A ring $R$ is said to be {\it pure-semi simple} if
every $R$-module is a direct sum of finitely generated
$R$-modules.
\end{definition}

\begin{theorem} The following are equivalent:\\
i) Every $R$-module is RD-pure injective.\\
ii) Every $R$-module is CP-injective.\\
iii) Every $R$-module is pure injective.\\
iv) $R$ is pure-semi simple.
\end{theorem}

{\bf Proof.} The implications $i)\Rightarrow ii)$ and
$ii)\Rightarrow iii)$ are clear.

Assume that iii) holds.  Then every pure exact sequence of
$R$-modules splits, and so it follows from [{\bf 8}] that every
$R$-module is a direct sum of finitely generated $R$-modules. Thus
iii) implies iv).

Now, we prove the implication $iv)\Rightarrow i)$. By [{\bf 7},
Theorem 4.3], $R$ is an Artinian principal ideal ring and every
$R$-module is a direct sum of cyclic $R$-modules. Hence, since
every ideal of $R$ is principal, it follows that every $R$-module
is a direct sum of modules of the form $R/Rr, r\in R$. From this
we can conclude that every RD-exact sequence splits. Therefore,
every $R$-module is RD-injective.  $\Box$

\section{Some examples}

Theorem 2.7 shows that there exists a ring $R$ such that
$\mathcal{C}_{RDR}\varsubsetneq\mathcal{C}_{CPR}\varsubsetneq
\mathcal{C}_{PR}
.$ In this section, we present some explicit
examples for these strict containments.

\begin{example} i) Let $\mathbb{Z}$ be the ring of integers and $p$ a
prime integer. Since every ideal of $\mathbb{Z}$ is principal, the
two notions of RD-injectivity and of CP-injectivity are coincide
for $\mathbb{Z}$-modules. Hence by [{\bf 3}, Theorem 3.6],
$D=\mathbb{Z}/p\mathbb{Z}$ is an RD-injective
$\mathbb{Z}$-module, while it is not an injective $\mathbb{Z}$-module.\\
ii) By [{\bf 1}, Example 1], there are an Artinian local ring $R$
and an $R$-algebra $S$ containing $R$, such that the inclusion map
$R\hookrightarrow S$ is cyclically pure, but it is not pure. It is
known that every Artinian $R$-module is pure injective (see e.g.
[{\bf 11}, Corollary 4.2]). Hence $R$ is a pure injective
$R$-module. But $R$ is not CP-injective, because otherwise by
[{\bf 3}, Theorem 3.4],  the inclusion map $R\hookrightarrow S$
splits.
\end{example}

\begin{lemma} Let $R$ be a domain, B a torsion free $R$-module
and $0\lo K\hookrightarrow B\lo M\lo 0$ an exact sequence of
$R$-modules. The following are equivalent:\\
i) $M$ is torsion-free.\\
ii) The inclusion map $K\hookrightarrow B$ is RD-pure.
\end{lemma}

{\bf Proof.} It is easy to see that an $R$-module $L$ is torsion
free if and only if $\Tor_1^R(R/Rr,L)=0$ for all $r\in R$. Since B
is torsion-free for any $r\in R$, from the exact sequence $0\lo
K\hookrightarrow B\lo M\lo 0$,  we deduce the exact sequence
$$0\lo \Tor^R_1(R/Rr,M)\lo (R/Rr)\otimes_RK\lo (R/Rr)\otimes_RB\lo
(R/Rr)\otimes_RM\lo 0.$$ Therefore, the assertion follows by
Remark 2.2 i). $\Box$

\begin{lemma} Let $R$ be a domain and $D$ an RD-injective
$R$-module. Then $\Ext_R^1(M,D)=0$ for all torsion-free
$R$-modules $M$.
\end{lemma}

{\bf Proof.} Let $M$ be a torsion-free $R$-module. Consider an
exact sequence $0\lo K\overset{i} \hookrightarrow F\lo M\lo 0 ,$
in which $F$ is a free $R$-module. Then, by Lemma 3.2, the
inclusion map $i$ is RD-pure. Now, from the exact sequence $$0\lo
\Hom_R(M,D)\lo \Hom_R(F,D)\lo \Hom_R(K,D)\lo \Ext_R^1(M,D)\lo 0,$$
we deduce that $\Ext_R^1(M,D)=0$. Note that since D is
RD-injective, the map $\Hom_R(i,id_D)$ is surjective. $\Box$

\begin{example} Let $(R,\fm)$ be a local Noetherian domain with
$\dim R>1$. Since $R$ is not a Pr\"{u}fer domain, it turns out
that $R$ possesses an ideal $\fa$ which is not projective. Thus
$\Ext_R^1(\fa,R/\fm)\neq 0$, by [{\bf 2}, Proposition 1.3.1]. Now,
by [{\bf 3}, Theorem 3.6], $R/\fm$ is a CP-injective $R$-module,
while by Lemma 3.3, $R/\fm$ is not RD-injective.
\end{example}

The following example shows that the two notions of
$\mathcal{S}_2$-flatness and $\mathcal{S}_2$-pureness are not the
same.

\begin{example} Assume that $R$ is a Noetherian domain such that $\dim
R>1$. Hence $R$ is not Pr\"{u}fer, and so by Corollary 2.12, there
exists an absolutely CP-module $M$ which is not injective. So,
there is an ideal $\fa$ such that $\Ext_R^ 1(R/{\fa},M)\neq 0$.
Let $E$ denote the injective envelope of $M$. Then from the exact
sequence $0\lo M\hookrightarrow E\overset{\pi}\lo E/M\lo 0 \ (*),$
we deduce the following exact sequence $$0\rightarrow
\Hom_R(R/{\fa},M)\rightarrow \Hom_R(R/{\fa},E)\rightarrow
\Hom_R(R/{\fa},E/M)\rightarrow \Ext_R^ 1(R/{\fa},M) \rightarrow
0.$$ Hence the map $\Hom_R(id_{R/{\fa}},\pi)$ is not surjective.
Thus $(*)$ is an $\mathcal{S}_2$-flat sequence which is not
$\mathcal{S}_2$-pure.
\end{example}


\end{document}